\documentclass[12pt]{amsart}

\usepackage{amssymb}

\usepackage{graphicx}
\usepackage{hyperref}

\newtheorem{theorem}{Theorem}[section]
\newtheorem{lemma}[theorem]{Lemma}

\theoremstyle{definition}

\newtheorem{proposition}[theorem]{Proposition}
\newtheorem{corollary}[theorem]{Corollary}

\newtheorem{remark}{Remark}

\theoremstyle{remark}

\numberwithin{equation}{section}
\usepackage{hyperref}

\begin{document}
\baselineskip=17pt
\title[Subset Containing $q^{th}$ Power Residue]{Almost No Finite Subset of Integers Contains a $q^{th}$ Power Modulo Almost Every Prime}


\author{Bhawesh Mishra}
\address{Department of Mathematical Sciences, The University of Memphis}
\curraddr{}
\email{mishra.bhawesh@memphis.edu}
\urladdr{}
\thanks{}


\date{\today}
\dedicatory{}

\commby{}

\begin{abstract}
Let $q$ be a prime. We give an elementary proof of the fact that for any $k\in\mathbb{N}$, the proportion of $k$-element subsets of $\mathbb{Z}$ that contain a $q^{th}$ power modulo almost every prime, is zero. This result holds regardless of whether the proportion is measured additively or multiplicatively. More specifically, the number of $k$-element subsets of $[-N, N]\cap\mathbb{Z}$ that contain a $q^{th}$ power modulo almost every prime is no larger than $a_{q,k} N^{k-(1-\frac{1}{q})}$, for some positive constant $a_{q,k}$. Furthermore, the number of $k$-element subsets of $\{\pm p_{1}^{e_{1}} p_{2}^{e_{2}} \cdots p_{N}^{e_{N}} : 0 \leq e_{1}, e_{2}, \ldots, e_{N}\leq N\}$ that contain a $q^{th}$ power modulo almost every prime is no larger than $m_{q,k} \frac{N^{Nk}}{q^{N}}$ for some positive constant $m_{q,k}$.
\end{abstract}

\maketitle

\section{Introduction}
We start our discussion with the notion of density, which, when defined, measures the proportion that a given $A\subseteq\mathbb{Z}$ occupies in the set of integers. 
\subsection{Additive and Multiplicative Density of Subsets}
The \textit{additive density} of a subset $A$ of $\mathbb{Z}$ is defined to be \[d(A) = \lim_{N\rightarrow\infty} \frac{|A \cap \{-N, -N+1, \ldots, -1, 0, 1, \ldots, N-1, N\}|}{2N+1},\] if the limit exists. Since $\mathbb{Z}$ is generated by $1$ as an additive group and $[-N, N]\cap\mathbb{Z}$ is the symmetric set containing the first $N$ multiples of its generator (including $0$), the quantity $d(A)$ is aptly called the additive density of $A$. Given a subset $A$ of integers and $k \in\mathbb{Z}$, define $A+k = \{a+k : a\in A\}$. If $d(A)$ exists, then for every $k\in\mathbb{Z}$, $d(A+k)$ also exists and is equal $d(A)$. In other words, the additive density is shift invariant. 

On the other hand, one can also define multiplicative density in integers, using the fact that $\mathbb{N}$, under multiplication, is a semi-group generated by the set of all the primes. Let $\mathbb{P} = \{p_{1} < p_{2} < \ldots < p_{i} < p_{i+1} < \ldots \}$ be the set of primes ordered by magnitude. We define the \textit{multiplicative density} of a subset $A$ of $\mathbb{Z}$ to be the quantity \[ d^{\times}(A) := \lim_{N\rightarrow\infty} \frac{|A \cap \{ \pm p_{1}^{e_{1}} p_{2}^{e_{2}} \cdots p_{N}^{e_{N}} : 0 \leq e_{1}, e_{2}, \ldots, e_{N} \leq N \}|}{2(N+1)^{N}},\]
if the limit exists. 

Just as the additive density of a subset (when it exists) is shift invariant, the multiplicative density of a subset (if it exists) is invariant to scaling by a non-zero integer. More specifically, if $d^{\times}(A)$ exists for $A \subseteq \mathbb{Z}$, then for every $k \in\mathbb{Z}\setminus\{0\}$, $d^{\times}(kA)$ exists and is equal to $d^{\times}(A)$. Here the set $kA$ is defined to be $\{m\in\mathbb{Z} : \frac{m}{k}\in A\}$. 

There are subsets of integers with additive density strictly between zero and one that have multiplicative density equal to $1$. For instance, $d(4\mathbb{Z}) = \frac{1}{4}$ because of shift-invariance of the additive-density and \[\mathbb{Z} = 4\mathbb{Z} \cup (4\mathbb{Z}+1) \cup (4\mathbb{Z}+2) \cup (4\mathbb{Z}+3).\] However, $1 = d^{\times}(\mathbb{Z}) = d^{\times}(4\mathbb{Z})$ due to the invariance of multiplicative density under scaling. On the other hand, \[d^{\times}\big( \bigcup_{i=1}^{3} (4\mathbb{Z}+i)\big) = d^{\times}\big(\mathbb{Z}\setminus4\mathbb{Z}\big) = 1 - d^{\times}(4\mathbb{Z}) = 0\] gives example of a set that has additive density $\frac{3}{4}$ but multiplicative density zero. In general, density can be defined in any amenable (semi)group, the class of which is quite rich.
For instance, the sub-class of amenable groups includes the class of all solvable groups, and is closed under directed unions, extensions and subgroups (see \cite[Ch. 10, Theorem 10.4]{Wagon}). 

For a subset $S$ of $\mathbb{P}$, the \textit{relative density} of $S$ (in $\mathbb{P}$) is defined to be the quantity \[ d_{\mathbb{P}}(S) := \lim_{N\rightarrow\infty} \frac{|S \cap \{1, 2, \ldots, N\}|}{|\mathbb{P} \cap \{1, 2, \ldots, N\}|}, \text{ if the limit exists.} \] We will say that a subset $A$ of integers contains a $q^{th}$ power modulo almost every prime if the set $\{p \in\mathbb{P} : A \text{ contains a } q^{th} \text{ power modulo } p \}$ has relative density $1$ in $\mathbb{P}$.

\subsection{Motivation}\label{motivation}
Unless otherwise stated, $q$ will denote a prime in this article. It is well known that if an integer $a$ is $q^{th}$ power modulo almost every prime, then $a$ is a perfect $q^{th}$ power in integers. Since there are at most $2\sqrt[q]{N}+1$ perfect $q^{th}$ powers in $\{-N, \ldots, -1, 0, 1, \ldots, N\}$, the set of integers that are $q^{th}$ power modulo almost every prime has additive density no more than \[\lim_{N\rightarrow\infty} \frac{2\sqrt[q]{N}+1}{2N+1} = \lim_{N\rightarrow\infty} \frac{1}{N^{1-\frac{1}{q}}} = 0.\] 

On the other hand, one can also show that the multiplicative density of perfect $q^{th}$ powers in the integers is zero. For instance, there are at most $\frac{N}{q}+1$ multiples of $q$ in the set $\{0, 1, \ldots, N\}$. Therefore, the set $\{\pm p_{1}^{e_{1}}p_{2}^{e_{2}}\cdots p_{N}^{e_{N}}: 0 \leq e_{1}, e_{2}, \ldots, e_{N}\leq N\}$ contains at most $2(\frac{N}{q}+1)^{N}$ perfect $q^{th}$ powers. Therefore, the multiplicative density of set of integers that are $q^{th}$ power modulo almost every prime, is at most
\begin{equation*}
        \lim_{N\rightarrow\infty} \frac{2(\frac{N}{q}+1)^{N}}{2(N+1)^{N}} = \lim_{N\rightarrow\infty} \frac{1}{q^{N}} = 0.
\end{equation*}
In other words, the proportion of $1$-element subset of $\mathbb{Z}$ that contain a $q^{th}$ power modulo almost every prime is zero, in both additive and multiplicative sense. 

It is proved in \cite[Corollary 1]{MishraFFA}  that only way a subset of integers with cardinality $\leq q$ can contain a $q^{th}$ power modulo almost every prime is when it contains a perfect $q^{th}$ power (see Corollary \ref{FFAcorollary} below). Therefore, using a similar argument, one can show that for $k \leq q$, the collection of $k$-element subsets of integers that contain a $q^{th}$ power modulo almost every prime has both additive and multiplicative subset density zero (see Lemmas \ref{qthpoweradditive} and \ref{qthpowermultiplicative}). 

However, when $k \geq q+1$, there are infinitely many $k$-element subsets of integers that contain a $q^{th}$ power modulo almost every prime (see Remark \ref{remark}). Therefore, it is not clear whether the proportion of $k$-element subsets that contain a $q^{th}$ power modulo almost every prime, is zero, when $k \geq q+1$. We show that this is in fact the case for every $k \in\mathbb{N}$, and give an elementary proof. 

Note that a $k$-element subset $\{a_{i}\}_{i=1}^{k}$ of integers contains a $q^{th}$ power modulo almost every prime if and only if the Galois group of the splitting field, over $\mathbb{Q}$, of the polynomial \[ f_{k}(x) = (x^{q} - a_{1}) (x^{q} - a_{2}) \cdots (x^{q} - a_{k})\] is equal to the union of conjugates of the stabilizers of the roots of its irreducible factors (see \cite[Theorem 1 and 2]{BerBil}). Such polynomials are also called \textit{exceptional} in the recent literature \cite{EKM}. Thus, the main result in this article is equivalent to the fact that for any $k\in\mathbb{N}$, the proportion of polynomials in $\mathbb{Z}[x]$ of the form \[ (x^{q} - a_{1}) (x^{q} - a_{2}) \cdots (x^{q} - a_{k})\] that are exceptional, is zero. Before we can state our main result, we require a generalization of additive and multiplicative densities, that allows us to measure the proportion of a given collection of subsets. 

\subsection{Additive and Multiplicative Subset Density}
Let $[N]$ denote the symmetric interval $\{-N, \ldots, -1, 0, 1, \ldots, N\}$ and for any $k \leq 2N+1$, let $\mathfrak{P}_{k}\big([N]\big)$ denote the collection of all $k$-element subsets of $[N]$. Let $\mathfrak{P}_{k}\big(\mathbb{Z}\big)$ denote the collection of every $k$-element subset of $\mathbb{Z}$. Then, the \textit{additive subset density} of a sub-collection $\mathcal{S}\subseteq\mathfrak{P}_{k}\big(\mathbb{Z}\big)$ is defined to be \[D(\mathcal{S}) := \lim_{N\rightarrow\infty} \frac{|\mathcal{S} \cap \mathfrak{P}_{k}\big([N]\big)|}{{2N+1 \choose k}} \text{, if the limit exists}.\] It is easily seen that the additive density of $A\subseteq\mathbb{Z}$ (if it exists) is equal to the additive subset density of the collection of $1$-element subsets $\big\{\{a\}\big\}_{a \in A}$. Therefore, the additive subset density of a collection of subsets can be viewed as a generalization of additive density. On the other hand, let \[ [N]^{\times} := \big\{ \pm p_{1}^{e_{1}}p_{2}^{e_{2}}\cdots p_{N}^{e_{N}} : 0 \leq e_{1}, e_{2}, \ldots, e_{N} \leq N \big\}. \] For any $k \leq 2(N+1)^{N}$, let $\mathfrak{P}_{k}\big([N]^{\times}\big)$ denote the collection of all $k$-element subsets of $[N]^{\times}$. The \textit{multiplicative subset density} of a sub-collection $\mathcal{T}\subset\mathfrak{P}_{k}(\mathbb{Z})$ is defined to be
\[D^{\times}\big(\mathcal{T}\big):= \lim_{N\rightarrow\infty} \frac{|\mathcal{T} \cap \mathfrak{P}_{k}\big([N]^{\times}\big)|}{{2(N+1)^{N} \choose k}}, \text{ if the limit exists.} \]

We will also repeatedly use that if $f$ and $g$ are two functions from $\mathbb{N}$ to $(0, \infty)$ such that $\lim_{N\rightarrow\infty} \frac{g(N)}{f(N)} = 0$, then for every real number $A$, there exists some constant $C > 0$ such that $f(N) + A \cdot g(N) \leq C\cdot f(N)$. Our main result is as follows:

\begin{theorem}
Let $q$ be a prime, $k$ be a natural number and $\mathcal{S}_{q,k}$ be the collection of $k$-element subsets of $\mathbb{Z}$ that contain a $q^{th}$ power modulo almost every prime. Then, the following hold:
\begin{enumerate}
    \item The additive subset density of $\mathcal{S}_{q,k}$ is zero. More specifically, the number of $k$-element subsets of $\{-N, \ldots, -1, 0, 1, \ldots, N\}$ that contain a $q^{th}$ power modulo almost every prime is no larger than $a_{q,k} N^{k-(1-\frac{1}{q})}$, for some positive constant $a_{q,k}$ depending upon $q$ and $k$. \vspace{2mm}

    \item The multiplicative subset density of $\mathcal{S}_{q,k}$ is zero. More specifically, the number of $k$-element subsets of $\{\pm p_{1}^{e_{1}} p_{2}^{e_{2}} \cdots p_{N}^{e_{N}} : 0 \leq e_{1}, e_{2}, \ldots, e_{N}\leq N\}$ that contain a $q^{th}$ power modulo almost every prime is no larger than $m_{q,k} \frac{N^{Nk}}{q^{N}}$, for some positive constant $m_{q,k}$ depending upon $q$ and $k$.
\end{enumerate}
\label{MainResult}
\end{theorem} 
This article has three more sections. In Section \ref{preliminary}, we will collect some results regarding sets that contain $q^{th}$ power modulo almost every prime, and prove some elementary lemmas. Section \ref{sectionsquare} contains the proof of Theorem \ref{MainResult} for $q = 2$, and Section \ref{sectionoddq} contains the same for odd $q$.

\section{Some Preliminary Results}\label{preliminary}
\subsection{Sets that contain a $q^{th}$ power modulo almost every prime}
The following classical result was first obtained by Fried in \cite{Fr} and was later rediscovered by Filaseta and Richman in \cite{FiRi}. 
\begin{proposition}
Let $A=\{a_{j}\}_{j=1}^{\ell}$ be a finite set of non-zero integers. Then, the following are equivalent:
\begin{enumerate}
    \item For almost every prime $p$, the set $A$ contains a square modulo $p$. 

    \item There exists a subset $S$ of $\{1, 2, \ldots, \ell\}$ with odd cardinality such that $\prod_{j \in S} a_{j}$ is a perfect square. 
\end{enumerate}
\label{Fried}
\end{proposition} 
We will need some definitions before stating a result akin to Proposition \ref{Fried} for odd primes. Let $q$ be an odd prime and $|a|$ denote the absolute value of an integer $a$. Given a positive integer $n > 1$ with unique prime factorization $\prod_{i=1}^{m} p_{i}^{a_{i}}$, we define the $q$-free part of $n$ to be \[ \text{rad}_{q}(n) := \prod_{i=1}^{m} p_{i}^{a_{i}(\text{mod } q)}.\] Let $\mathbb{F}_{q}$ denote the finite field with $q$-elements and $V$ be a vector space over $\mathbb{F}_{q}$. A collection of subspaces $\{W_{j}\}_{j\in J}$ of $V$ is said to be a \textit{linear covering} of $V$ if $\bigcup_{j\in J} W_{j} = V$. We will say that a linear covering $\{W_{j}\}_{j\in J}$ of $\mathbb{F}_{q}^{r}$ is a covering by hyperplanes if for every $j \in J$, there exists $\nu_{1j}, \nu_{2j}, \ldots, \nu_{rj}\in\mathbb{F}_{q}$, not all zero, such that \[W_{j} := \Big\{ \big(x_{i}\big)_{i=1}^{r} : \sum_{i=1}^{r} \nu_{ij}x_{i} = 0 \Big\}.\]
\subsubsection{Linear Hyperplanes Associated with a Finite Subset of Integers}\label{hyperplanedefn}
Let $q$ be an odd prime and $B = \{b_{j}\}_{j=1}^{\ell}$ be a finite subset of integers that does not contain a perfect $q^{th}$ power. Let $p_{1}, p_{2}, \ldots, p_{r}$ be all the primes that divide any element of the set $\{\text{rad}_{q}\big(|b_{j}|\big) : 1\leq j\leq\ell\}$. For every $1 \leq i \leq r$ and every $1 \leq j\leq \ell$, let $p_{i}^{\nu_{ij}}\mid\mid \text{rad}_{q}(|b_{j}|)$ for some $\nu_{ij} \geq 0$. We will say that the set 
\[\Big\{ \big\{ \big(x_{i})_{i=1}^{r}: \sum_{i=1}^{k} \nu_{ij}x_{i} = 0 \big\} : 1 \leq j \leq \ell \Big\}\] is the set of hyperplanes in $\mathbb{F}_{q}^{r}$ associated with the set of integers $B = \{b_{j}\}_{j=1}^{\ell}$. We will use the following proposition. 
\begin{proposition}
Let $q$ be an odd prime and $B=\{b_{j}\}_{j=1}^{\ell}$ be a finite subset of integers that does not contain a perfect $q^{th}$ power. Then, the following conditions are equivalent:
\begin{enumerate}
    \item For almost all primes $p$, the set $B$ contains a $q^{th}$ power modulo $p$. 

    \item For every $(c_{j})_{j=1}^{\ell}\in\mathbb{Z}^{\ell}$, there exists two subsets $A, B$ of $\{1, 2, \ldots, \ell\}$ with $|A| \not\equiv |B| \hspace{1mm} (\text{mod } q)$ such that $\prod_{j\in A} b_{j} = g^{q} \prod_{j\in B} b_{j}$ for some non-zero integer $g$. 

    \item The set of hyperplanes \[\Big\{ \big\{ \big(x_{i})_{i=1}^{r}: \sum_{i=1}^{k} \nu_{ij}x_{i} = 0 \big\} : 1 \leq j \leq \ell \Big\}\] associated with $B$, forms a linear covering of $\mathbb{F}_{q}^{r}$. 
\end{enumerate}
\label{FFAresult}
\end{proposition}
The equivalence between the conditions $(1)$ and $(2)$ in Proposition \ref{FFAresult} is a consequence of \cite[Theorem 10.2]{Sk}, which itself was a refinement of a much more general result obtained in \cite{SS} by Schinzel and Ska\l{}ba. The equivalence between $(1)$ and $(3)$ appears in \cite[Theorem 1]{MishraFFA}. A consequence of the equivalence of $(1)$ and $(3)$ is the following corollary (see \cite[Corollary 1]{MishraFFA} for a proof).
\begin{corollary} 
Let $q$ be an odd prime and $k \in\mathbb{N}$ such that $k \leq q$. A $k$-element subset, $B$, of integers contains a $q^{th}$ power modulo every prime if and only if $B$ contains a perfect $q^{th}$ power. \label{FFAcorollary}
\end{corollary}

\begin{remark}\label{remark}
Let $k = q+1$ and $p_{1}, p_{2}\in\mathbb{P}\setminus\{q\}$ be two distinct primes. One can easily see that \[ \mathcal{H} := \Big\{ (x_{1}, x_{2})\in\mathbb{F}_{q}^{2} :  x_{1} = 0 \Big\} \bigcup_{j=0}^{q-1} \Big\{ (x_{1}, x_{2})\in\mathbb{F}_{q}^{2} : jx_{1}+x_{2} = 0 \Big\} \] forms a linear covering of $\mathbb{F}_{q}^{2}$, because any $(x_{1}, x_{2})\in\mathbb{F}_{q}^{2}$ satisfies the equation of at least one hyperplane in $\mathcal{H}$. Therefore, Proposition \ref{FFAresult} implies that the set \[ \mathcal{B} := \big\{ p_{1}, p_{1}^{j} p_{2} : j \in\{0, 1, \ldots, q-1\} \big\} \] contains a $q^{th}$ power modulo almost every prime. Therefore, there exists infinitely many subsets of integers of cardinality $\geq q+1$ that contain a $q^{th}$ power modulo almost every prime. 
\end{remark}

Before we proceed with the proof of Theorem \ref{MainResult}, we will state two elementary lemmas, which will bound the number of $k$-element subsets, respectively of $[N]$ and $[N]^{\times}$, that contain at least one perfect $q^{th}$ power.

\begin{lemma}
Let $q$ be a prime and $k, N\in\mathbb{N}$. Then, the number of $k$-element subset of $[N]$ that contain at least one perfect $q^{th}$ power is no more than $a^{\prime}_{q,k} N^{k-(1-\frac{1}{q})}$, for some positive constant $a^{\prime}_{q,k}$ depending upon $q$ and $k$.  
\label{qthpoweradditive}
\end{lemma}

\begin{proof}
When $q$ is odd, there are at most $2\sqrt[q]{N}+1$ many perfect $q^{th}$ powers in the set $[N]$. For every $1 \leq i \leq k$, there are ${2\sqrt[q]{N}+1 \choose i}$ many ways of choosing $i$ perfect $q^{th}$ powers from $[N]$ and at most ${2N \choose k-i}$ ways of choosing the remaining $k-i$ elements that are not perfect $q^{th}$ powers. Therefore, the total number of $k$-element subsets of $[N]$ that contain at least one perfect $q^{th}$ power is at most $\sum_{i=1}^{k} {2\sqrt[q]{N}+1 \choose i}{2N \choose k-i}$. 

For $q = 2$, there are at most $\sqrt{N}+1$ perfect squares in $[N]$. Therefore, the total number of $k$-element subsets of $[N]$ that contain at least one perfect square is at most $\sum_{i=1}^{k} {\sqrt{N}+1 \choose i}{2N \choose k-i}$. The term with the highest power of $N$ in the above sums is
\[\begin{cases}
    \frac{2^{k-1}}{(k-1)!} N^{k-1/2} \text{ when } q = 2, \\
    \frac{2^{k}}{(k-1)!} N^{k-(1-\frac{1}{q})} \text{ when } q \text{ is odd}.
\end{cases}\]
Therefore, the number of $k$-element subsets that contain at least one perfect $q^{th}$ power is no larger than $a^{\prime}_{q,k} N^{k-(1-\frac{1}{q})}$, for some positive constant $a^{\prime}_{q,k}$ depending on both $q$ and $k$. 
\end{proof}

\begin{lemma}
Let $q$ be a prime and $k, N\in\mathbb{N}$. Then, the number of $k$-element subsets of $[N]^{\times}$ that contain at least one perfect $q^{th}$ power is no larger than $m^{\prime}_{q,k} \cdot \frac{N^{Nk}}{q^{N}}$, for some positive constant $m^{\prime}_{q,k}$ depending upon $q$ and $k$.  
\label{qthpowermultiplicative}
\end{lemma}

\begin{proof}
When $q$ is odd, there are at most $\frac{N}{q}+1$ many multiples of $q$ in $\{0, 1, 2, \ldots, N\}$. So, there are at most $2 \big(\frac{N}{q} + 1)^{N}$ perfect $q^{th}$ powers, in $[N]^{\times}$. There are at most ${2 \big(\frac{N}{q} + 1)^{N}+ 1 \choose i}$ ways of choosing $i$ perfect $q^{th}$ powers from $[N]^{\times}$, and at most ${2(N+1)^{N} \choose k-i}$ ways of choosing the rest of the elements that are not perfect $q^{th}$ powers. Hence, there are at most \[\sum_{i=1}^{k} {2 \big(\frac{N}{q} + 1)^{N}+ 1 \choose i} {2(N+1)^{N} \choose k-i}\] ways of choosing a $k$-element subset of $[N]^{\times}$ that contain at least one perfect $q^{th}$ power. The term with the highest power of $N$ in the ${2 \big(\frac{N}{q} + 1)^{N}+ 1 \choose i} {2(N+1)^{N} \choose k-i}$ is equal to \[\frac{2^{i} N^{Ni}}{q^{Ni} \cdot i!} \cdot \frac{2^{k-i} N^{N(k-i)}}{(k-i)!} = \frac{2^{k}N^{Nk}}{q^{Ni} \cdot i! (k-i)!}.\] Therefore, the term with the highest power of $N$ in $ \sum_{i=1}^{k} {2 \big(\frac{N}{q} + 1)^{N}+ 1 \choose i} {2(N+1)^{N} \choose k-i}$ is 
\begin{multline*}
\\ 2^{k} N^{Nk} \sum_{i=1}^{k} \frac{1}{{q^{Ni} i! (k-i)!}} \leq 2^{k}N^{Nk} \sum_{i=1}^{k} \frac{1}{q^{N}}  = \frac{2^{k} \cdot k \cdot N^{Nk}}{q^{N}}.\\
\end{multline*}
Therefore, the number of $k$-element subset of $[N]^{\times}$ that contain at least one perfect $q^{th}$ power is $m^{\prime}_{q,k} \cdot \frac{N^{Nk}}{q^{N}}$, for some positive constant $m^{\prime}_{q,k}$ depending on $q$ and $k$. 
 
For $q = 2$, there are at most $\frac{N}{2}+1$ multiples of $2$ in $\{0, 1, \ldots, N\}$ and hence at most $(\frac{N}{2}+1)^{N}$ perfect squares in $[N]^{\times}$. Therefore by exactly analogous argument, we obtain that there are at most \[\sum_{i=1}^{k} {{\big(\frac{N}{2} + 1)^{N} \choose i}} {2(N+1)^{N} \choose k-i}\] many $k$-element subsets of $[N]^{\times}$ that contain at least one perfect $q^{th}$ power. The term with the highest power of $N$ in ${\big(\frac{N}{2} + 1)^{N} \choose i} {2(N+1)^{N} \choose k-i}$ is \[\frac{N^{Ni}}{2^{Ni}\cdot i!} \cdot \frac{2^{k-i} N^{N(k-i)}}{(k-i)!} = \frac{2^{k} \cdot N^{Nk}}{2^{Ni+i} \cdot i! (k-i)!}, \] and hence the term with the highest power of $N$ in $\sum_{i=1}^{k} {\big(\frac{N}{2} + 1)^{N}+ 1 \choose i} {2(N+1)^{N} \choose k-i}$ is \[ 2^{k} \cdot N^{Nk} \sum_{i=1}^{k} \frac{1}{2^{i(N+1)} \cdot i! (k-i)!} \leq 2^{k} \cdot N^{Nk} \sum_{i=1}^{k} \frac{1}{2^{N+1}} \leq \frac{k2^{k-1} N^{Nk}}{2^{N}},\] which shows that the number of $k$-element subsets of $[N]^{\times}$ that contain at least one perfect square is no more than $m^{\prime}_{2,k} \cdot \frac{N^{Nk}}{2^{N}}$ for some constant depending on $k$. 
\end{proof}

Now we are ready to prove our main result. 
\section{Proof of Theorem \ref{MainResult} for $q = 2$}\label{sectionsquare}
According to the Proposition \ref{Fried}, if $B \in\mathcal{S}_{2,k}$ and $B$ does not contain a perfect square, then there exists $A \subseteq B$ of odd cardinality $\geq 3$ such that $\prod_{a \in A} a$ is a perfect square. Let $A = \{a_{\mu}\}_{\mu=1}^{r}$ be such a subset of $B$, then the condition $(2)$ in Proposition \ref{Fried} implies that for every $j\in\{1, 2, \ldots, r\}$, we have \[\text{rad}_{2}(a_{j}) = \text{rad}_{2}\Big(\prod_{\mu=1, \mu\neq j}^{r} a_{j}\Big).\] In other words, fixing any $r-1$ elements of $A$ determines the square-free part of the last remaining element. Now, we will divide this section into two further subsections, to deal with additive and multiplicative subset density separately. 
\subsection{Proof that $\mathcal{S}_{2,k}$ has additive subset density zero.}  
We have already shown in Lemma \ref{qthpoweradditive} that the number of $k$-element subsets of $[N]$ that contain a perfect square is no more than $a^{\prime}_{2,k} \cdot N^{k-\frac{1}{2}}$. Now, we will find a bound for the number of $k$-element subsets of $[N]$ that do not contain a perfect square, yet end up containing a square modulo almost every prime. 
 
There are at most ${2N \choose r-1}$ ways of choosing $r-1$ integers from $[N]$, that are not perfect squares. For any such choice of $r-1$ elements of $A$, one has a unique choice for the square-free part of the remaining $r^{th}$ element. Therefore, there are at most $\sqrt{N}+1$ choices for the last remaining element of $A$. Furthermore, there are exactly ${2N+1-r \choose k-r}$ ways of choosing the remaining $(k-r)$ elements of $B$. Therefore, there are at most \[{2N \choose r-1}\cdot {2N+1-r \choose k-r} \cdot (\sqrt{N}+1)\] many $k$-element subsets $B$ of $[N]$ such that the following is satisfied:
\begin{itemize}
    \item $B$ contains a square modulo almost every prime but does not contain a perfect square.  

    \item $r$ is the lowest odd number $\leq k$ such that $B$ contains a subset $A$ of cardinality $r$ and $\prod_{a\in A} a$ is a perfect square. 
\end{itemize}

The total number of $k$-element subsets of $[N]$ that contain a square modulo almost every prime, but do not contain a perfect square is thus at most
\begin{equation}
\sum_{r=3, r \text{ odd}}^{k} {2N \choose r-1}\cdot {2N+1-r \choose k-r} \cdot (\sqrt{N}+1).
\label{nosquare}
\end{equation}
Note that the term with the highest power of $N$ in ${2N \choose r-1}\cdot {2N+1-r \choose k-r} \cdot (\sqrt{N}+1)$ is \[ \frac{2^{r-1} N^{r-1}}{(r-1)!} \frac{2^{k-r} N^{k-r}}{(k-r)!} \sqrt{N} =  \frac{2^{k-1}}{(r-1)!(k-r)!}N^{k-\frac{1}{2}}\] and hence the term with the highest power of $N$ that appears in 
\[\sum_{r=3, r \text{ odd}}^{k} {2N+1 \choose r-1}\cdot {2N+1-r \choose k-r} \cdot (\sqrt{N}+1)\] is 
\[\gamma_{2,k} N^{k-\frac{1}{2}}, \text{ where } \gamma_{2,k}:= 2^{k-1} \sum_{r=3, r \text{ odd }}^{k} \frac{1}{(r-1)!(k-r)!}.  \]

Thus, the number of $k$-element subsets of $[N]$ that do not contain a perfect square, but do contain a square modulo almost every prime is at most $a^{\prime\prime}_{2,k}$ $N^{k-\frac{1}{2}}$, for some positive constant $a^{\prime\prime}_{2,k}$ depending on $k$. This fact, along with Lemma \ref{qthpoweradditive} for $q = 2$ implies that the number of $k$-element subsets of $[N]$ that contain a square modulo almost every prime is no larger than \[a^{\prime\prime}_{2,k} \cdot N^{k-\frac{1}{2}} + a^{\prime}_{2,k} \cdot N^{k-\frac{1}{2}} = a_{2,k} N^{k-\frac{1}{2}}, \text{ where } a_{2,k} = a^{\prime}_{2,k} + a^{\prime\prime}_{2,k},\] which establishes the upper-bound in $(1)$ of Theorem \ref{MainResult} when $q = 2$. Therefore,
\begin{multline*}
 \\ \limsup_{N\rightarrow\infty} \frac{|\mathcal{S}_{2,k} \cap \mathfrak{P}_{k}\big([N]\big)|}{{2N+1 \choose k}} \leq \limsup_{N\rightarrow\infty} \frac{a_{2,k} \cdot N^{k-\frac{1}{2}}}{\frac{2^{k}\cdot N^{k}}{k!}} = \frac{k! \cdot a_{2,k}}{2^{k}} \lim_{N\rightarrow\infty} \frac{1}{N^{\frac{1}{2}}} = 0, \\ 
\end{multline*}
which establishes that the additive subset density of $\mathcal{S}_{2,k}$ is zero.

\subsection{Proof that $\mathcal{S}_{2,k}$ has multiplicative subset density zero.}
We have already shown in Lemma \ref{qthpowermultiplicative} that the number of $k$-element subsets of $[N]^{\times}$ that contain a perfect square is no more than $m^{\prime}_{2,k} \cdot \frac{N^{Nk}}{2^{N}}$. Now, we will find a bound for the number of $k$-element subsets of $[N]^{\times}$ that do not contain a perfect square, yet end up containing a square modulo almost every prime.

Note that there are at most ${2(N+1)^{N} \choose r-1}$ ways of choosing $r-1$ integers from $[N]^{\times}$ that are not perfect squares and for any such choice of $r-1$ elements of $A$, we have a unique choice for the square-free part of the last remaining element, which can be multiplied by a non-zero perfect square. Therefore, there are $(\frac{N}{2}+1)^{N}$ many choices for the last remaining element, once the $r-1$ elements are chosen. Furthermore, there are exactly ${2(N+1)^{N}-r \choose k-r}$ ways of choosing the remaining $(k-r)$ elements of $B$. 

Therefore, there are at most \[{2(N+1)^{N} \choose r-1}\cdot {2(N+1)^{N}-r \choose k-r} \cdot (\frac{N}{2}+1)^{N}\] many $k$-element subsets $B$ of $[N]^{\times}$ such that the following is satisfied:
\begin{itemize}
    \item $B$ contains a square modulo almost every prime but does not contain a perfect square.  

    \item $r$ is the lowest odd number $\leq k$ such that $B$ contains a subset $A$ of cardinality $r$ and $\prod_{a\in A} a$ is a perfect square. 
\end{itemize}
The total number of $k$-element subsets of $[N]^{\times}$ that contain a square modulo almost every prime, but does not contain a perfect square is obtained again by summing over all odd numbers $r$ such that $3 \leq r \leq k$, i.e,
\begin{equation}
\sum_{r=3, r \text{ odd}}^{k} {2(N+1)^{N} \choose r-1}\cdot {2(N+1)^{N}-r \choose k-r} \cdot (\frac{N}{2}+1)^{N}.
\label{nosquaremult}
\end{equation}

Note that the term with highest power of $N$ in \[ {2(N+1)^{N} \choose r-1}\cdot {2(N+1)^{N}-r \choose k-r} \cdot (\frac{N}{2}+1)^{N}\] is
\begin{equation*}
    \frac{2^{r-1}\cdot N^{N(r-1)}}{(r-1)!} \cdot \frac{2^{k-r}\cdot N^{N(k-r)}}{(k-r)!} \cdot \frac{N^{N}}{2^{N}} = \frac{N^{Nk}}{2^{N}} \cdot \frac{2^{k-1}}{(r-1)!(k-r)!}.
\end{equation*} 
Therefore, the term with highest power of $N$ in the sum \[\sum_{r=3, r \text{ odd}}^{k} {2N^{N+1} \choose r-1}\cdot {2N^{N+1}-r \choose k-r} \cdot (\frac{N}{2}+1)^{N}\] is equal to 
\[\eta_{2,k} \cdot \frac{N^{Nk}}{2^{N}}, \text{ where } \eta_{2,k}:= \sum_{r=3, r \text{ odd}}^{k} \frac{2^{k-1}}{(r-1)!(k-r)!}.\]

In other words, the number of $k$-element subsets of $[N]^{\times}$ that do not contain a perfect square, but do contain a square modulo almost every prime is at most $m^{\prime\prime}_{2,k} \frac{N^{Nk}}{2^{N}}$, for some positive constant $m^{\prime\prime}_{2,k}$ depending upon $k$. This fact, along with Lemma \ref{qthpowermultiplicative} for $q = 2$ implies that the number of $k$-element subsets of $[N]$ that contain a square modulo almost every prime is no larger than \[m^{\prime\prime}_{2,k} \frac{N^{Nk}}{2^{N}} + m^{\prime}_{2,k} \frac{N^{Nk}}{2^{N}} = m_{2,k} \frac{N^{Nk}}{2^{N}}, \text{ where } m_{2,k} = m^{\prime\prime}_{2,k}+m^{\prime}_{2,k}\] which establishes the upper bound in $(2)$ of Theorem \ref{MainResult} when $q=2$. Therefore, 
\begin{multline*}
\\ \limsup_{N\rightarrow\infty} \frac{|\mathcal{S}_{2,k} \cap \mathfrak{P}_{k}([N]^{\times})|}{{2(N+1)^{N} \choose k}} \leq \limsup_{N\rightarrow\infty} \frac{m_{2,k} \frac{N^{Nk}}{2^{N}}}{\frac{2^{k} N^{Nk}}{k!}} = \frac{k! \cdot m_{2,k}}{2^{k}} \lim_{N\rightarrow\infty} \frac{1}{2^{N}} = 0,\\
\end{multline*}
which proves that the multiplicative subset density of $\mathcal{S}_{2,k}$ is zero too.

\section{Proof of Theorem \ref{MainResult} when $q$ is an odd prime}\label{sectionoddq}
We will again divide the proof into two subsections each of which will deal with additive and multiplicative subset density separately.

\subsection{Proof that $D(\mathcal{S}_{q,k}) = 0$}\label{subsectionadditiveoddq}
We have already shown in Lemma \ref{qthpoweradditive} that the number of $k$-element subsets of $[N]$ that contain a perfect $q^{th}$ power is no more than $a^{\prime}_{q,k} \cdot N^{k-(1-\frac{1}{q})}$. Now, we will find a bound for the number of $k$-element subsets $B$ of $[N]$ that do not contain a perfect $q^{th}$ power, yet end up containing a $q^{th}$ power modulo almost every prime. 

If $D = \{d_{j}\}_{j=1}^{k}$ is a $k$-element subset of $[N]$, not containing a perfect $q^{th}$ power, then using the condition $(2)$ in Proposition \ref{FFAresult}, with $(c_{j})_{j=1}^{k} = (1, 1, \ldots, 1)$, we obtain that there exists two subsets $A, B$ of $\{1, 2, \ldots, k\}$ with $|A| \not\equiv |B| \hspace{1mm} (\text{mod } q)$ and $\prod_{j\in A} d_{j} = g^{q} \prod_{j \in B} d_{j}$ for some non-zero integer $g$. 

More specifically, if the elements of $A$ and $|B|-1$ elements of $B$ are fixed, then there is a unique $q$-free choice for the last element of $B$. We can also assume that $A \cap B = \emptyset$, because if not, one can always factor out the common $d_{j}$ from the both side of $\prod_{j\in A} d_{j} = g^{q} \prod_{j \in B} d_{j}$, for every $j\in A \cap B$. 

Assume that $|A| = r$ and $|B| = s$ with $r \not\equiv s \hspace{1mm} (\text{mod } q)$. Note that there are at most ${2N+1 \choose r+s-1}$ ways of choosing $r$ elements of $A$ and $(s-1)$ elements of $B$. Then, the $q$-free part of last remaining element of $B$ is fixed, which we can modify by at most $2\sqrt[q]{N}+1$ perfect $q^{th}$ powers. Furthermore, the remaining $(k-r-s)$ elements of $D$ can be chosen to be any non-zero elements. Therefore, there are at most ${2N+1 \choose k-r-s}$ choices for remaining elements of $D$. Hence, there are at most \[{2N+1 \choose r+s-1} {2N+1 \choose k-r-s} (2\sqrt[q]{N}+1)\] ways of obtaining a $k$-element subset $D=\{d_{j}\}_{j=1}^{k}$ of $[N]$ such that:
\begin{enumerate}
    \item $|A| = r$, $|B| = s$, $r \not\equiv s \hspace{1mm} (\text{mod } q)$, $A \cap B = \emptyset$ and

    \item $\prod_{j\in A} d_{j} = g^{q} \prod_{j\in B} d_{j}$.
\end{enumerate}
Therefore, the total number of $k$-element subsets $D$ of $[N]$ that do not contain a perfect $q^{th}$ power but contains $q^{th}$ power modulo almost every prime is at most
\begin{equation}
    \sum_{r=s=1, r\not\equiv s \hspace{1mm}(\text{mod } q)}^{k} {2N+1 \choose r+s-1} {2N+1 \choose k-r-s} (2\sqrt[q]{N}+1). \label{oddqestimate1}
\end{equation}
The term with the highest power of $N$ in ${2N+1 \choose r+s-1} {2N+1 \choose k-r-s} (2\sqrt[q]{N}+1)$ is equal to 
\begin{multline*}
   \\ \frac{2^{r+s-1} N^{r+s-1}}{(r+s-1)!} \cdot \frac{2^{k-r-s} N^{k-r-s}}{(k-r-s)!} 2\sqrt[q]{N} = \frac{2^{k}}{(r+s-1)!(k-r-s)!} N^{k-(1-\frac{1}{q})}\\
\end{multline*}
Therefore, the term with the highest power of $N$ in the sum \eqref{oddqestimate1} is equal to \[ \gamma_{q,k} N^{k-(1-\frac{1}{q})}, \text{ where } \gamma_{q,k} := \sum_{r=s=1, r\not\equiv s \hspace{1mm}(\text{mod } q)}^{k} \frac{2^{k}}{(r+s-1)!(k-r-s)!}. \] Therefore, the total number of $k$-element subsets $D$ of $[N]$ not containing a perfect $q^{th}$ power but containing a $q^{th}$ power modulo every prime is at most $a^{\prime\prime}_{q,k} N^{k-(1-\frac{1}{q})}$, for some positive constant $a^{\prime\prime}_{q,k}$ depending upon $q$ and $k$. 

Now, using Lemma \ref{qthpoweradditive} for odd $q$, we therefore get hat the total number of $k$-element subsets of $[N]$ that contain a $q^{th}$ power modulo almost every prime is no larger than 
\[a^{\prime\prime}_{q,k} N^{k-(1-\frac{1}{q})} + a^{\prime}_{q,k} N^{k-(1-\frac{1}{q})} \leq a_{q,k} N^{k-(1-\frac{1}{q})},\] where $a_{q,k} = a^{\prime\prime}_{q,k} + a^{\prime}_{q,k}$. Thus,
\begin{multline*}
 \limsup_{N\rightarrow\infty} \frac{|\{\mathcal{S}_{q,k} \cap \mathfrak{P}_{k}([N])\}|}{{2N+1 \choose k}} \leq \limsup_{N\rightarrow\infty} \frac{a_{q,k} N^{k-(1-\frac{1}{q})}}{\frac{2^{k} N^{k}}{k!}} = \frac{a_{q,k} \cdot k!}{2^{k}} \lim_{N\rightarrow\infty} \frac{1}{N^{1-\frac{1}{q}}} = 0, 
\end{multline*}
and hence, the additive subset density of $\mathcal{S}_{q,k}$ is zero. 

\subsection{Proof that $D^{\times}(\mathcal{S}_{q,k}) = 0$}
Exactly as in subsection \ref{subsectionadditiveoddq}, the discussion analogous to that before \eqref{oddqestimate1} gives that the number of $k$ element subsets of $[N]^{\times}$ that do not contain a perfect $q^{th}$ power but contains a $q^{th}$ power modulo almost every prime, is no larger than 
\begin{equation}
    \sum_{r=s=1, r\not\equiv s \hspace{1mm}(\text{mod } q)}^{k} {2(N+1)^{N} \choose r+s-1} {2(N+1)^{N} \choose k-r-s} (\frac{N}{q}+1)^{N}. \label{oddqestimate2}
\end{equation}
Note that \eqref{oddqestimate2} is obtained from \eqref{oddqestimate1} by replacing $|[N]| = 2N+1$ by $|[N]^{\times}| = 2(N+1)^{N}$, and the maximum number of perfect $q^{th}$ powers in $[N]$, i.e., $2\sqrt[q]{N}+1$ by the maximum number of perfect $q^{th}$ powers in $[N]^{\times}$, i.e., $(\frac{N}{q}+1)^{N}$. 

Note that the term with the highest power of $N$ in \[{2(N+1)^{N} \choose r+s-1} {2(N+1)^{N} \choose k-r-s} (\frac{N}{q}+1)^{N} \text{ is } \frac{2^{k}}{(r+s-1)!(k-r-s)!} \frac{N^{Nk}}{q^{N}}.\] So the term with the highest power of $N$ in \eqref{oddqestimate2} is equal to \[ \eta_{q,k} \frac{N^{Nk}}{q^{N}}, \text{ where } \eta_{q,k} := \sum_{r=s=1, r\not\equiv s \hspace{1mm}(\text{mod } q)}^{k} \frac{2^{k}}{(r+s-1)!(k-r-s)!}. \]

Hence, the number of $k$-element subsets of $[N]^{\times}$ that do not contain a perfect $q^{th}$ power but contain a $q^{th}$ power is at most $m^{\prime\prime}_{q,k} \frac{N^{Nk}}{q^{N}}$, for some positive constant $m^{\prime\prime}_{q,k}$ depending on both $q$ and $k$. Hence, using Lemma \ref{qthpowermultiplicative} for odd $q$ gives that
\[|\mathcal{S}_{q,k} \cap \mathfrak{P}_{k}([N])^{\times}| \leq m^{\prime\prime}_{q,k} \frac{N^{Nk}}{q^{N}} + m^{\prime}_{q,k} \frac{N^{Nk}}{q^{N}} \leq m_{q,k} \frac{N^{Nk}}{q^{N}}, \] where $m_{q,k} = m^{\prime\prime}_{q,k} + m^{\prime}_{q,k}$. Thus,
\begin{multline*}
\limsup_{N\rightarrow\infty} \frac{|\mathcal{S}_{q,k} \cap \mathfrak{P}_{k}([N]^{\times})|}{{2(N+1)^{N} \choose k}} \leq \limsup_{N\rightarrow\infty} \frac{m_{q,k} \frac{N^{Nk}}{q^{N}}}{\frac{2^{k}\cdot N^{Nk}}{k!}} = \frac{k! \cdot m_{q,k}}{2^{k}} \lim_{N\rightarrow\infty} \frac{1}{q^{N}} = 0,
\end{multline*}
which completes the proof of $(2)$ in Theorem \ref{MainResult}.

\bibliographystyle{plain}
\bibliography{FA}
\end{document}